\input amstex
\magnification=1200
\documentstyle{amsppt}
\NoRunningHeads
\NoBlackBoxes
\define\CDb{\Cal D}
\define\Ff{\frak f}
\define\Fg{\frak g}
\define\Fh{\frak h}
\define\CU{\Cal U}
\define\sla{\operatorname{\frak s\frak l}}
\define\sLtwo{\sla(2,\Bbb C)}
\define\RW{\operatorname{\Cal R\Cal W}}
\define\cRW{\widehat{\RW}}
\define\RWtwo{\RW(\sLtwo)}
\define\cRWtwos{\cRW'(\sLtwo)}
\define\cRWtwo{\cRW(\sLtwo)}
\define\SRW{S(\cRWtwo,\sLtwo)}
\define\ZRW{Z(\cRWtwo,\sLtwo)}
\define\Res{\operatorname{Res}}
\hyphenation{Zhe-lo-ben-ko}
\document\qquad\qquad\qquad\qquad\qquad\qquad\qquad\qquad\qquad\qquad
$\boxed{\boxed{\aligned
&\text{\eightpoint"Thalassa Aitheria" Reports}\\
&\text{\eightpoint RCMPI-96/06 (October 1996)}\endaligned}}$\newline
\ \newline
\ \newline
\topmatter
\title Topics in hidden symmetries. II
\endtitle
\author Denis V. Juriev
\endauthor
\affil\eightpoint\rm "Thalassa Aitheria" Research Center for Mathematical
Physics and Informatics,\linebreak
ul.Miklukho-Maklaya 20-180, Moscow 117437, Russia.\linebreak
E-mail: denis\@juriev.msk.ru
\endaffil
\date E--version: q-alg/yymmxxx\enddate
\abstract This short paper being devoted to some aspects of the inverse
problem of the representation theory briefly treats the interrelations
between the author's approach to the setting free of hidden symmetries and
the researches of D.P.Zhelobenko in the generalized Mickelsson algebras
and related topics.
\endabstract
\endtopmatter
This short paper being a continuation of the first part [1] is a collection of
examples illustrating the general ideology presented in the review [2]. The
examples should emphasize the interrelations between the material of the
previous papers (see f.e.[3] and [1:\S2]) on the setting free of hidden
symmetries and the researches of D.P.Zhelobenko on the generalized Mickelsson
algebras [4,5] and related objects of the representation theory, so the
material may be considered also as an illustrative commentary to the book [5].

\head 1. Topic Four: The generalized Mickelsson S--algebra $\SRW$ and the
Zhelobenko Z--algebra $\ZRW$ of the central extension $\cRWtwo$ of the
Racah--Wigner algebra $\RWtwo$ for the Lie algebra $\sLtwo$
\endhead

This topic is devoted to the generalized Mickelsson S--algebra and the
Zhelobenko Z--algebra of the central extension $\cRWtwo$ of the Racah--Wigner
algebra $\RWtwo$, which was introduced in [6:\S2.2;3] (see also
[2:\S\S1.2,1.3;7:\S1.2]).

The algebra $\cRWtwo$ is generated by the nine elements $l_i$ ($i=-1,0,1$),
$w_i$ ($i=-2,-1,0,1,2$) and $\varrho$, the least is a central element.
The commutation relations have the form
$$\allowdisplaybreaks\align
[l_i,l_j]&=(i-j)l_{i+j},\\
[l_i,w_j]&=(2i-j)w_{i+j},\\
[w_{-1},w_{-2}]&=2l_{-1}\circ w_{-2},\\
[w_0,w_{-2}]&=\tfrac43(2l_{-1}\circ w_{-1}+l_0\circ w_{-2}),\\
[w_0,w_{-1}]&=\tfrac16(-l_1\circ w_{_2}+10l_0\circ w_{-1}+3l_{-1}\circ w_0-
3\varrho l_{-1}),\\
[w_1,w_{-2}]&=l_1\circ w_{-2}+2l_0\circ w_{-1}+3l_{-1}\circ w_0+\varrho l_{-1},\\
[w_1,w_{-1}]&=\tfrac12(l_1\circ w_{-1}+6l_0\circ w_0+l_{-1}\circ w_1-\varrho l_0),\\
[w_2,w_{-2}]&=4(l_1\circ w_{-1}+l_{_1}\circ w_1),\\
[w_2,w_{-1}]&=2l_{-1}\circ w_2+2l_0\circ w_1+3l_1\circ w_0+\varrho l_1,\\
[w_1,w_0]&=\tfrac16(-l_{-1}w_2+10l_0\circ w_1+3l_1\circ w_0-3\varrho l_1),\\
[w_2,w_0]&=\tfrac43(2l_1\circ w_1+l_0\circ w_2,\\
[w_2,w_1]&=2l_1\circ w_2,
\endalign$$
where $X\circ Y=\frac12(XY+YX)$. Note that one may use the right or the left
ordering for the products of $l_i$ and $w_j$ in the r.h.s. instead of their
Weyl ordering (it should be marked also that the formulas in [6:\S2.2] contain
misprints in two commutators, which were corrected in [3]). The hypothesis of
[6:\S2.2] suggests that the spectrum of the degenerate Verma modules of the
finitely generated associative algebra $\cRWtwo$ of polynomial growth is deeply
related to the Kac spectrum for the Virasoro algebra [8,9] and possibly may be
used for the description of some models ("quasi--minimal") of two--dimensional
field theory with broken infinite--dimensional Virasoro symmetries (cf.[10]).

\remark{Remark 1} The associative algebra $\cRWtwo$ is neither contragredient
nor even finitely generated locally triangular algebra in sense of
D.P.Zhelobenko [5:App.A] (see also [11]).
\endremark

The generators $l_i$ form the Lie algebra $\sLtwo$ so $\cRWtwo$ is an
envelopping algebra for $\sLtwo$, moreover, it admits an embedding of
$\CU(\sLtwo)$. The pair $(\cRWtwo,\sLtwo)$ obeys the conditions of
D.P.Zhelobenko [4,5:\S7.3], so one may construct the generalized Mickelsson
S--algebra $\SRW$ and the Zhelobenko Z--algebra $\ZRW$ in lines of [4,5].

Let us denote the span of $l_0$ by $\Fh$, $\CU(\Fh)$ and $\CDb(\Fh)$ are the
universal envelopping algebra of $\Fg$ and its algebra of quotients [12],
$\cRWtwos=\cRWtwo\otimes_{\CU(\Fh)}\CDb(\Fh)$. Then, the generalized Mickelsson
S--algebra $\ZRW$ is the algebra over $\CDb(\Fh)$ (in sense of D.P.Zhelobenko)
generated by the six generators $v_i$ ($i=-2,-1,0,1,2$) and $r$, the least is
central one. The generators are the operators $pw_i$ and $p\varrho$ in $\cRWtwos$
reduced to the quotient $\cRWtwos/I_+$, where $I_+$ is the left ideal generated
by $l_1$, and $p$ is the extremal projector for $\sLtwo$ [5:\S\S3.2,3.3].
The extremal projector has the form:
$$p=\sum_{n=0}^{\infty}\frac{(-1)^n}{n!(2l_0-2)\ldots(2l_0-n-1)}l_{-1}^nl_1^n,
\text{\ so that \ } l_1p=pl_{-1}=0,$$
(this formula for $p$ slightly differs from one of D.P.Zhelobenko in view of
the different choices of basises in $\sLtwo$).

\proclaim{Theorem 1} The Zhelobenko Z--algebra $\ZRW$ is defined
by the quadratic (over $\CDb(\Fh)$) relations between the generators $v_i$
($i=-2,-1,\mathbreak 0,1,2$) and $r$ (the least is central one). The relations
between $v_i$ have the form:
$$\allowdisplaybreaks\align
v_{-1}v_{-2}&=(1-\tfrac2{\eta-2})v_{-2}v_{-1},\\
v_0v_{-2}&=(1-\tfrac6{(\eta-2)(2\eta-3)})v_{-2}v_0-
\tfrac4{\eta-1}v_{-1}^2+\tfrac{4\eta}3v_{-2},\\
v_0v_{-1}&=(1-\tfrac3{\eta-1})v_{-1}v_0+
\tfrac{2(\eta-3)}{(\eta-2)(2\eta-3)}v_{-2}v_1+\tfrac{5\eta-2}3v_{-1},\\
v_1v_{-2}&=(1-\tfrac6{\eta(\eta-1)})v_{-2}v_1+
\tfrac{12(\eta-2)}{\eta(2\eta-1)}v_{-1}v_0+
\tfrac{2(\eta-1)(\eta-2)}{\eta}v_{-1},\\
v_1v_{-1}&=(1-\tfrac9{(\eta-1)(2\eta-1)})v_{-1}v_1+
\tfrac{\eta(2\eta-5)}{2(\eta-1)(\eta-2)(2\eta-3)}v_{-2}v_2-
\tfrac9{2\eta}+\tfrac{3(2\eta+1)}2v_0+\tfrac{\eta}2r,\\
v_2v_{-2}&=(1-\alpha)v_{-2}v_2-
\tfrac{2(\eta-3)(\eta+2)}{\eta(\eta-1)(\eta+1)}v_{-1}v_1+
\tfrac{36}{(\eta+1)(2\eta+1)}v_0^2-\tfrac{12\eta}{\eta+1}v_0+
\tfrac{\eta(\eta-3)}{\eta+1}r,\\
v_2v_{-1}&=(1-\tfrac6{\eta(\eta+1)})v_{-1}v_2+
\tfrac{12(\eta-1)}{(\eta+1)(2\eta+1)}v_0v_1+
\tfrac{2\eta(\eta-1)}{\eta+1}v_1,\\
v_1v_0&=(1-\tfrac3{\eta})v_0v_1-
\tfrac{2(\eta-2)}{(\eta-1)(2\eta-1)}v_{-1}v_2+\tfrac{5\eta+3}3v_1,\\
v_2v_0&=(1-\tfrac6{\eta(2\eta+1)})v_0v_2-
\tfrac4{\eta+1}v_1^2+\tfrac{4(\eta+2)}3v_2,\\
v_2v_1&=(1-\tfrac2{\eta+1})v_1v_2,
\endalign$$
where $\eta=l_0\in\CDb(\Fh)$ (note that $[\eta,v_i]=-iv_i$) and
$\alpha=\tfrac4{(\eta+1)(\eta-2)}(1-\tfrac9{2(2\eta-1)(2\eta-3)})$.
\endproclaim

The proof is just in lines of the proof of Theorem 4.2.4 in [5:\S4.2] with
slight but evident modifications. The concrete computations are made in the
standard way (see [5:\S4.2]) here.

\remark{Remark 2} The described construction of the generalized Mickelsson
S--algebra $\SRW$ and the Zhelobenko Z--algebra $\ZRW$ admits a
superanalogue. The role of the Racah--Wigner algebra $\RWtwo$ for the Lie
algebra $\sLtwo$ is played by the Racah--Wigner algebra $\RW(\sla(2|1,\Bbb C))$
for the Lie superalgebra $\sla(2|1,\Bbb C))$ [6:\S2.3].
\endremark

\head 2. Topic Five: The generalized Mickelsson S--algebras $S(\mho(\Fg,\pi),\Fg)$
and the Zhelobenko Z--algebras $Z(\mho(\Fg,\pi),\Fg)$ of the mho--algebras
$\mho(\Fg,\pi)$ for the complex semisimple Lie algebras $\Fg$
\endhead

This topic is devoted to the generalized Mickelsson S--algebras and the
Zhelobenko Z--algebras of the mho-algebras $\mho(\Fg,\pi)$ for the complex
semisimple Lie algebras $\Fg$ and their finite--dimensional representations,
which were introduced and discussed in [1:\S2] (see also [2:\S1.4]).

\definition{Definition}

{\bf A.} Let $\Fg$ be a Lie algebra and $\pi$ be its (irreducible)
representation. {\it Mho--algebra\/} $\mho(\Fg,\pi)$ is an associative algebra
such that (1) $\CU(\Fg)$ is a subalgebra of $\mho(\Fg,\pi)$ and, hence, $\Fg$
naturally acts in $\mho(\Fg,\pi)$, (2) there is defined a $\Fg$--equivariant
embedding of $\pi$ into $\mho(\Fg,\pi)$, so $\pi$ may be considered as
a subspace of $\mho(\Fg,\pi)$, (3) the $\Fg$--equivariant embedding of $\pi$
into $\mho(\Fg,\pi)$ is extended to a $\Fg$--equivariant embedding of
$S^{\cdot}(\pi)$ into $\mho(\Fg,\pi)$, defined by the Weyl symmetrization,
and, therefore, $S^{\cdot}(\pi)$ may be considered as a subspace of
$\mho(\Fg,\pi)$; (4) $\Fg$--modules $\mho(\Fg,\pi)$ and $S^{\cdot}(\Fg)\otimes
S^{\cdot}(\pi)$ are isomorphic, here the isomorphism of subalgebra $\CU(\Fg)$
of the algebra $\mho(\Fg,\pi)$ and $S^{\cdot}(\Fg)$ as $\Fg$--modules is
used; (5) in an arbitrary basis $w_l$ in $\pi$ the commutator of two elements
of the basis in the algebra $\mho(\Fg,\pi)$ may be represented in the form
$[w_i,w_j]=f^k_{ij}w_k$, where the "noncommutative structural functions"
$f^k_{ij}$ are the elements of the algebra $\CU(\Fg)$.

{\bf B.} Let $\Fg$ be a Lie algebra and $\pi$ be its (irreducible)
representation. {\it Affine mho--algebra\/} $\hat\mho(\Fg,\pi)$ is
an associative algebra such that the conditions (1)--(4) above hold and
(5') in an arbitrary basis $w_l$ in $\pi$ the commutator of two elements of the
basis in the algebra $\hat\mho(\Fg,\pi)$ may be represented in the form
$[w_i,w_j]=f^k_{ij}w_k+g_{ij}$, where the "noncommutative structural
functions" $f^k_{ij}$ and $g_{ij}$ are the elements of the algebra $\CU(\Fg)$.
\enddefinition

Notations $\mho(\Fg,\pi)$ and $\hat\mho(\Fg,\pi)$ emphasize an analogy
between mho--algebras and affine mho--algebras and the universal envelopping
algebras. The examples of the mho--algebras for $\Fg=\sLtwo$ were considered
in [1:\S2] (see also [2:\S1.4]). For instance, the natural semi--direct product
$\CU(\sLtwo)\ltimes\CU(\sLtwo)$ and the Racah--Wigner algebra $\RWtwo$ are
mho--algebras, whereas $\CU(\sLtwo\oplus\sLtwo)$ with the diagonal embedding
of $\CU(\sLtwo)$ and $\cRWtwo$ are affine mho--algebras. The mho-algebra
$\mho(\sLtwo,\pi_3)$ ($\pi_3$ is the seven--dimensional representation of
$\sLtwo$) was considered in [1:\S2] (see also [2:\S1.4]).

\proclaim{Theorem 2} Let $\Fg$ be a complex semisimple Lie algebra,
and $\mho(\Fg,\pi)$ be an arbitrary mho-algebra over $\Fg$, then the pair
$(\mho(\Fg,\pi),\Fg)$ obeys the conditions of D.P.Zhe\-lo\-ben\-ko [4;5:\S7.3],
the generalized Mickelsson S--algebra $S(\mho(\Fg,\pi),\Fg)$, the
Zhe\-lo\-ben\-ko Z--algebra $Z(\mho(\Fg,\pi),\Fg)$ may be constructed as in [4,5],
the least is generated by the elements $px$, $x\in\pi$, $p$ is the extremal
projector for $\Fg$, and it is a quadratic algebra over $\CDb(\Fh)$ in sense
of D.P.Zhelobenko ($\Fh$ is the Cartan subalgebra of $\Fg$).

This statement holds true for affine mho--algebras $\hat\mho(\Fg,\pi)$ also.
\endproclaim

The simplest way to restore the proof of this theorem for the reader is to
calculate explicitely the quadratic relations in the Zhelobenko Z--algebra
$Z(\mho(\sLtwo,\pi_3),\mathbreak\sLtwo)$ (the commutation relations in
$\mho(\sLtwo,\pi_3)$ were written by the author in [1:\S2], see also [2:\S1.4]).

\remark{Remark 3} One may consider the reductive Lie algebra $\Fg$ instead
of the semisimple one. In this case the mentioned above conditions of
D.P.Zhelobenko should be checked to mantain the statement of the theorem 2.
\endremark

\remark{Remark 4} One may also consider the generalized Mickelsson
S--algebras $S(\mho(\Fg,\pi),\Ff)$ and the Zhelobenko Z--algebras
$Z(\mho(\Fg,\pi),\Ff)$, where $\Ff$ is a reductive
algebra reductively embed into $\Fg$ such that $\Res^{\Fg}_{\Ff}(\pi)$ is
completely reducible ($\Res$ denotes the restriction functor here).
\endremark
\newpage

\Refs
\roster
\item" [1]" Juriev D., Topics in hidden symmetries. I.: E-print hep-th/9405050
(1994).
\item" [2]" Juriev D., An excursus into the inverse problem of representation
theory [in Russian]: Report RCMPI-95/04 (August 1995) [e-version:
mp\_arc/96-477 (1996)].
\item" [3]" Juriev D., Setting hidden symmetries free by the noncommutative
Veronese mapping. J.Math.Phys. 35(9) (1994) 5021-5024.
\item" [4]" Zhelobenko D.P., Extremal projectors and generalized Mickelsson
algebras over reductive Lie algebras [in Russian]. Izvestiya AN SSSR.
Ser.matem. 52(4) (1988) 758-773.
\item" [5]" Zhelobenko D.P., Representations of the reductive algebras.
Moscow, Nauka, 1994.
\item" [6]" Juriev D., Complex projective geometry and quantum projective field
theory [in Russian]. Teor.Matem.Fiz. 101(3) (1994) 331-348 [English transl.:
Theor.Math.Phys. 101 (1994) 1387-1403].
\item" [7]" Juriev D., Infinite dimensional geometry and quantum field theory
of strings. III. Infinite dimensional W--geometry of a second quantized
free string. J.Geom.Phys. 16 (1995) 275-300 [e-version: hep-th/9401026].
\item" [8]" Kac V.G., Infinite dimensional Lie algebras. Cambridge, Cambridge
Univ. Press, 1990.
\item" [9]" Feigin B.L., Fuchs D.B., Representations of the Virasoro algebra.
In "Representations of infinite dimensional Lie algebras". Gordon and Breach,
1991.
\item"[10]" Fradkin E.S., Palchik M.Ya., Exactly solvable models of
conformal--invariant quantum field theory in D--dimensional space. J.Geom.Phys.
5(4) (1988) 601-629 [reprinted in "Geometry and physics. Essays in honor of
I.M.Gelfand", Eds.S.Gindikin and I.M.Singer, Bologna, Pitagora Editrice and
Amsterdam, Elsevier Sci. Publ., 1991].
\item"[11]" Zhelobenko D.P., Contragredient algebras. J.Group Theory Phys.
1(1) (1993) 201-233.
\item"[12]" Dixmier J., Alg\`ebres enveloppantes. Paris, Villars, 1973.
\endroster
\endRefs
\enddocument